\documentclass[10pt]{amsart}

\theoremstyle{remark}
\newtheorem*{rem}{Remark}

\usepackage{url}
\usepackage[colorlinks,linkcolor=blue,citecolor=blue]{hyperref}

\title[A simple computation of $\zeta(2k)$]{A simple computation of $\zeta(2k)$ by using Bernoulli polynomials and a telescoping series}

\author[\'O.~Ciaurri]{\'Oscar Ciaurri}
\address{%
        Depar\-ta\-men\-to de Mate\-m\'a\-ti\-cas y Com\-pu\-ta\-ci\'on,
        Uni\-ver\-si\-dad de La Rioja,
        26004 Lo\-gro\-\~no, Spain}
\email{oscar.ciaurri@unirioja.es}

\author[L.~M.~Navas]{Luis M. Navas}
\address{%
        Depar\-ta\-men\-to de Mate\-m\'a\-ti\-cas,
        Uni\-ver\-si\-dad de Salamanca,
        37008 Sala\-manca, Spain}
\email{navas@usal.es}

\author[F.~J.~Ruiz]{Francisco J. Ruiz}
\address{%
        Depar\-ta\-men\-to de Mate\-m\'a\-ti\-cas,
        Uni\-ver\-si\-dad de Zaragoza,
        50009 Zara\-goza, Spain}
\email{fjruiz@unizar.es}

\author[J.~L.~Varona]{Juan L. Varona}
\address{%
        Depar\-ta\-men\-to de Mate\-m\'a\-ti\-cas y Com\-pu\-ta\-ci\'on,
        Uni\-ver\-si\-dad de La Rioja,
        26004 Lo\-gro\-\~no, Spain}
\email{jvarona@unirioja.es}
\urladdr{http://www.unirioja.es/cu/jvarona/}

\thanks{The research of the authors is supported by grant MTM2012-36732-C03-02 of the~DGI}



\date{\textit{Amer.\ Math.\ Monthly}.}

\begin{document}

\begin{abstract}
We present a new simple proof of Euler's formulas for $\zeta(2k)$, where $k = 1,2,3,\dots$ 
The computation is done using only the defining properties of the Bernoulli polynomials 
and summing a telescoping series, and the same method also yields integral formulas 
for $\zeta(2k+1)$.
\end{abstract}

\maketitle

\section{Introduction}

In the mathematical literature, one finds many ways of obtaining the formula
\begin{equation}
\label{eq:sum2k}
  \zeta(2k) := \sum_{n=1}^{\infty} \frac{1}{n^{2k}}
             = \frac{(-1)^{k-1}2^{2k-1}\pi^{2k}}{(2k)!} B_{2k},
             \qquad k = 1,2,3,\dots,
\end{equation}
where $B_k$ is the $k$th Bernoulli number, a result first published by Euler in 1740. For example, the recent paper~\cite{AmoDF} contains quite a complete list of references; among them, the articles \cite{Apo-MO, Os-MO, Ts-MO, Wi-MO} published in this \textsc{Monthly}.
The aim of this paper is to give a new proof of~\eqref{eq:sum2k} which is simple and elementary, in the sense that it involves only basic one variable Calculus, the Bernoulli polynomials, and a telescoping series. As a bonus, it also yields integral formulas for $\zeta(2k+1)$ and the harmonic numbers.

\subsection{The Bernoulli polynomials --- necessary facts}

For completeness, we begin by recalling the definition of the Bernoulli polynomials $B_{k}(t)$ and their basic properties. There are of course multiple approaches one can take (see~\cite{CDG}, which shows seven ways of defining these polynomials). A frequent starting point is the generating function
\begin{equation*}
  \frac{x e^{xt}}{e^x-1} = \sum_{k=0}^\infty B_k(t)\, \frac{x^k}{k!},
\end{equation*}
from which, by the uniqueness of power series expansions, one can quickly obtain many of their basic properties. Among these we single out
\begin{gather}
\label{eq:deriv}
  B_{0}(t) = 1,
  \quad
  B'_k(t) = k B_{k-1}(t),
  \qquad k\geq 1,
\end{gather}
which shows by induction that $B_{k}(t)$ is in fact a polynomial, and
\begin{equation}
\label{eq:integ}
  \int_{0}^{1} B_{k}(t) \, dt = 0, \qquad k \geq 1.
\end{equation}
Alternatively, one can instead use~\eqref{eq:deriv} and~\eqref{eq:integ} to define the polynomials $B_{k}(t)$ recursively.
In any case, one finds that the first few Bernoulli polynomials are
\[
  B_0(t) = 1,
  \qquad
  B_1(t) = t - \frac{1}{2},
  \qquad
  B_2(t) = t^2 - t + \frac{1}{6}.
\]
The Bernoulli \emph{numbers} are defined to be the values $B_k = B_k(0)$. From 
~\eqref{eq:deriv}, \eqref{eq:integ}, and the Fundamental Theorem of Calculus, one sees that $B_{k}(0) = B_{k}(1)$ for $k \geq 2$, and from the symmetry relation
\begin{equation*}
  B_k(1-t) = (-1)^k B_k(t),
  \qquad k\geq 0
\end{equation*}
(easily proved by induction on~$k$) one deduces
\begin{equation}
\label{eq:symBnum}
 \begin{aligned}
    B_{2k}(1) 
 &= B_{2k}(0) = B_{2k},
 \qquad
 k \geq 0,
 \\
    B_{2k+1}(1) 
 &= B_{2k+1}(0) = 0,
 \qquad
 k \geq 1.
 \end{aligned}
\end{equation}
%
Of course there are many other properties and relations satisfied by the Bernoulli polynomials, but those listed above are the only ones necessary for our goal.

\subsection{Outline of the proof}

Given these basic facts about Bernoulli polynomials, let us give a sketch of our proof of \eqref{eq:sum2k} (the details are in the next section). The integrals
\[
  I^{*}(k,m) := \int_0^{1} B_{2k}^{*}(t) \cos(m\pi t)\,dt
\]
where $B_{k}^{*}(t) = B_{k}(t) - B_{k}(0)$, are evaluated via a recurrence formula obtained
from integrating by parts twice. Solving the recurrence and summing over $m$ gives
\[
  \frac{(-1)^{k-1}(2k)!}{2^{2k} \pi^{2k}} \, \zeta(2k) = \sum_{m=1}^{\infty} I^{*}(k,2m).
\]
The key step in the proof is to apply the elementary trigonometric identities relating products of sines and cosines to sums to obtain a formula for $\cos(m x)$ which expresses the latter sum as a telescoping series of integrals. Another integration by parts justifies passage to the limit in these integral representations (showing that the general term tends to zero) and, together with \eqref{eq:integ}, yields
\[
  \sum_{m=1}^{\infty} I^{*}(k,2m) = \frac{1}{2}\,B_{2k},
\]
which gives Euler's formula.

\subsection{Other applications of the ideas of the proof}

The same technique, applied to the odd integer case, gives us the integral expression
\cite[formula 23.2.17]{AS}
\[
  \zeta(2k+1) = \frac{(-1)^{k-1}2^{2k}\pi^{2k+1}}{(2k+1)!}
  \int_0^{1} B_{2k+1}(t) \cot\Bigl(\frac{\pi t}{2}\Bigr)\,dt,
\]
with terms that mimic \eqref{eq:sum2k} except for the innocent-looking yet nonetheless completely mysterious integral, for which there is no known ``nice'' closed form expression (e.g., one which could determine the irrationality or even the transcendence of all values $\zeta(2k+1)$). 

In spite of this state of affairs, there is no lack of formulas for $\zeta(2k+1)$ in the mathematical literature. There is, for example, an intriguing parametric formula due to Ramanujan in which the Bernoulli numbers appear. For positive $\alpha,\beta$ with $\alpha\beta=\pi^2$ and $k$ any nonzero integer, we have
\begin{multline*}
  \alpha^{-k} \left( \tfrac12 \zeta(2k+1) + \sum_{n=1}^{\infty} \frac{n^{-2k-1}}{e^{2\alpha n}-1} \right) 
  = (-\beta)^{-k} \left( \tfrac12 \zeta(2k+1) + \sum_{n=1}^{\infty} \frac{n^{-2k-1}}{e^{2\beta n}-1} \right) \\
  - 2^{2k} \sum_{n=0}^{k+1} (-1)^n \,\frac{B_{2n}}{(2n)!} \frac{B_{2k+2-2n}}{(2k+2-2n)!}\, \alpha^{k+1-n} \beta^n,
\end{multline*}
see \cite[Entry 21 (i) on page 275]{Ber}; in the recent papers \cite{MSW, GMR}, this formula has been analyzed from the standpoint of transcendence.
The book \cite[Section 4.2]{SrCh} contains a large collection of other formulas for $\zeta(2k+1)$.

\section{Computation of $\zeta(2k)$}

\subsection{Some auxiliary integrals}

Consider the integrals
\[
  I(k,m) := \int_0^{1} B_{2k}(t) \cos(m\pi t)\,dt,
  \qquad
  k \geq 0,\ m \geq 1.
\]
An immediate computation shows that $I(0,m) = 0$ for $m \geq 1$. For $k \geq 1$, 
integrating by parts twice and applying~\eqref{eq:deriv}, we get
\begin{align*}
     I(k,m) 
  &=   \frac{1}{m\pi}\Bigl[B_{2k}(t) \sin(m \pi t)\Bigr]_{t=0}^{t=1}
     - \frac{2k}{m\pi} \int_0^{1} B_{2k-1}(t) \sin(m\pi t)\,dt
  \\
  &=   \frac{2k}{m^{2} \pi^{2}} \Bigl[B_{2k-1}(t) \cos(m \pi t) \Bigr]_{t=0}^{t=1}
     - \frac{2k(2k-1)}{m^{2} \pi^{2}} \,I(k-1,m),
\end{align*}
which gives us both the special case
\[
  I(1,m) = \int_0^{1} \left( t^2 - t + \frac{1}{6} \right) \cos(m\pi t)\,dt =
  \begin{cases}
  0, & m=1,3,5,\dots, \\
  \dfrac{2}{m^2 \pi^2}, & m=2,4,6,\dots
  \end{cases}
\]
and, by~\eqref{eq:symBnum}, the recurrence relation
\[
 I(k,m) = - \frac{2k(2k-1)}{m^{2} \pi^{2}} \,I(k-1,m),
 \qquad
 k \geq 2.
\]
From this recurrence one easily obtains the closed form
\begin{equation}
\label{eq:evalInt}
    I(k,m) = 
  \begin{cases}
  0, & m=1,3,5,\dots, \\
  \dfrac{(-1)^{k-1}(2k)!}{m^{2k} \pi^{2k}}, & m=2,4,6,\dots.
  \end{cases}
\end{equation}

Now, for reasons which are made clear below, consider $B^{*}_{k}(t) = B_{k}(t) - B_{k}(0) = B_{k}(t) - B_{k}$, i.e., the Bernoulli polynomial minus its constant term. The corresponding integral
\[
  I^{*}(k,m) := \int_0^{1} B_{2k}^{*}(t) \cos(m\pi t)\,dt
              = \int_0^{1} (B_{2k}(t)-B_{2k}) \cos(m\pi t)\,dt
\]
is equal to $I(k,m)$, because $\int_0^{1} \cos(m\pi t)\,dt = 0$ for $m>0$. For fixed $k \geq 1$, summing~\eqref{eq:evalInt} over $m$ yields
\begin{align*}
  \frac{(-1)^{k-1}(2k)!}{2^{2k} \pi^{2k}} \, \zeta(2k)
  &= \frac{(-1)^{k-1}(2k)!}{\pi^{2k}} \sum_{m=1}^{\infty} \frac{1}{(2m)^{2k}}
  = \sum_{m=1}^{\infty} I^{*}(k,2m) 
  = \sum_{m=1}^{\infty} I^{*}(k,m).
\end{align*}

\subsection{The telescoping trick}

We will need the elementary trigonometric identity
\begin{equation}
\label{eq:telescope}
  \cos(mx) = \frac{\sin(\frac{2m+1}{2}x)-\sin(\frac{2m-1}{2}x)}{2\sin(\frac{x}{2})}.
\end{equation}
With the introduction of~\eqref{eq:telescope}, we now have a telescoping series, yielding
\begin{multline*}
     \frac{(-1)^{k-1}(2k)!}{2^{2k} \pi^{2k}} \, \zeta(2k) 
     = \sum_{m=1}^{\infty} \int_0^{1} B^{*}_{2k}(t) \cos(m\pi t)\,dt
  \\
  \begin{aligned}
  &= \lim_{N \to \infty} \sum_{m=1}^{N} \Bigg( 
    \int_0^{1} B^{*}_{2k}(t) \frac{\sin(\frac{2m+1}{2}\pi t)}{2\sin(\frac{\pi t}{2})}\,dt
    - \int_0^{1} B^{*}_{2k}(t) \frac{\sin(\frac{2m-1}{2}\pi t)}{2\sin(\frac{\pi t}{2})}\,dt
  \Bigg) 
  \\
  &= \left(\lim_{N \to \infty} \int_0^{1} B^{*}_{2k}(t) 
    \frac{\sin(\frac{2N+1}{2}\pi t)}{2\sin(\frac{\pi t}{2})}\,dt\right)
  - \frac{1}{2} \int_0^{1} B^{*}_{2k}(t) \,dt.
  \end{aligned}
\end{multline*}
We observe that by~\eqref{eq:integ}, the value of the second term is
\[
    \frac12 \int_0^{1} B^{*}_{2k}(t) \,dt
  = \frac12 \int_0^{1} (B_{2k}(t)-B_{2k}) \,dt
  = - \frac{B_{2k}}{2}.
\]
Now, we show that the limit in the first term is $0$. Note that the function
\[
  f(t) = \frac{B^{*}_{2k}(t)}{2\sin(\frac{\pi t}{2})}, \qquad t \in (0,1],
\] 
extends by continuity to $t=0$ since $B^{*}_{2k}(0)=0$ (this is the reason for subtracting the constant term), and is differentiable on $[0,1]$ with a continuous derivative. Denoting $(2N + 1)\pi/2$ by~$R$, integrating by parts gives
\[
  \int_0^1 f(t) \sin(Rt)\,dt 
  = - \frac{\cos(R)}{R} f(1) + \frac{1}{R} f(0) + \int_0^1 f'(t) \frac{\cos(Rt)}{R} \,dt.
\]
The boundedness of $f'(t)$ shows that each term in the above sum approaches zero as $R \to \infty$, so that indeed the limit tends to~$0$.
Consequently,
\[
  \frac{(-1)^{k-1}(2k)!}{2^{2k}\pi^{2k}} \,\zeta(2k) = \frac{B_{2k}}{2},
\]
which, after rearranging terms, gives~\eqref{eq:sum2k}.

\section{What about $\zeta(2k+1)$?}

The same approach will yield a formula for $\zeta(2k+1)$, but the term which is subtracted when summing the telescoping series is an integral which, as far as anyone knows, cannot be evaluated in a simple closed form. We proceed in the same way as before, except this time, we consider the integrals
\[
 J(k,m) := \int_0^{1} B_{2k+1}(t) \sin(m\pi t)\,dt.
\]
Direct computation shows that
\[
  J(0,m) = \int_0^{1} \left( t - \frac{1}{2} \right) \sin(m\pi t)\,dt 
         = - \frac{1 + (-1)^{m}}{m \pi}
         = 
  \begin{cases}
  0,                 & m=1,3,5,\dots, 
                     \\
  -\dfrac{1}{m \pi}, & m=2,4,6,\dots.
  \end{cases}
\]
For $k \geq 1$, integrating by parts twice gives, using \eqref{eq:deriv} and \eqref{eq:symBnum}, the recurrence relation
\[
 J(k,m) = - \frac{(2k+1)(2k)}{m^2 \pi^{2}} \,J(k-1,m).
\]
From this recurrence we obtain the closed form
\begin{equation}
\label{eq:evalJ}
  J(k,m) = 
  \begin{cases}
  0, & m=1,3,5,\dots, \\
  \dfrac{(-1)^{k-1}(2k+1)!}{m^{2k+1} \pi^{2k+1}}, & m=2,4,6,\dots.
  \end{cases}
\end{equation}
Note that this time, since $B_{2k+1}(0) = B_{2k+1} = 0$ for $k \geq 1$, subtracting the constant term is not necessary except for $k = 0$, since $B_{1} = -1/2$, although, since $\zeta(1) = \infty$, this is irrelevant for now (but see the next Remark). Thus, by \eqref{eq:evalJ}, and using the trigonometric identity 
\begin{equation}
\label{eq:telescopeSin}
  \sin(mx) = -\frac{\cos(\frac{2m+1}{2}x)-\cos(\frac{2m-1}{2}x)}{2\sin(\frac{x}{2})},
\end{equation}
we obtain, for $k \geq 1$,
\begin{multline*}
     \frac{(-1)^{k-1}(2k+1)!}{2^{2k+1} \pi^{2k+1}} \,\zeta(2k+1)
  = \frac{(-1)^{k-1}(2k+1)!}{\pi^{2k+1}} \sum_{m=1}^{\infty} \frac{1}{(2m)^{2k+1}}
  = \sum_{m=1}^{\infty} J(k,m)
  \\
  = -\left(\lim_{N \to \infty} \int_0^{1} B_{2k+1}(t) 
  \frac{\cos(\frac{2N+1}{2}\pi t)}{2\sin(\frac{\pi t}{2})}\,dt\right)
  + \int_{0}^{1} B_{2k+1}(t) \frac{\cos(\frac{\pi t}{2})}{2\sin(\frac{\pi t}{2})}\,dt.
\end{multline*}
The limit is null, for the same reason as before. Consequently, we have proved that
\begin{equation}
\label{eq:zeta-impar}
  \zeta(2k+1) = \frac{(-1)^{k-1}2^{2k}\pi^{2k+1}}{(2k+1)!}
  \int_0^{1} B_{2k+1}(t) \cot\Bigl(\frac{\pi t}{2}\Bigr)\,dt,
  \qquad k  \geq 1.
\end{equation}
%
It would be nice to know if the integral in \eqref{eq:zeta-impar} has a closed form expression (other than in terms of $\zeta(2k+1)$ of course!) but at present this problem remains open.


\begin{rem}
In the case $k=0$, corresponding to the harmonic series $\sum_{m=1}^{\infty} 1/m = \infty$, we can still obtain information from the integrals $J(0,m)$ if we consider the partial sums of the telescoping series.
The trick is to just use \eqref{eq:telescopeSin} to sum $\sin(m \pi x)$ \emph{inside} the integral $J(0,m)$. 
This leads to the formula
\begin{equation}
\label{eq:harmonic}
  H_{M}
 = \pi \int_{0}^{1} \left(t - \frac{1}{2}\right) \frac{\cos\left(\frac{4M+1}{2} \pi t\right) 
  - \cos\left(\frac{\pi t}{2}\right)}{\sin\left(\frac{\pi t}{2}\right)} \, dt
\end{equation}
where $H_{M} = \sum_{m=1}^{M} {1}/{m}$ is the $M$th harmonic number.
\end{rem}

\section{Fourier Confidential}

Readers familiar with the basic theory of Fourier series will of course immediately recognize that the integrals $I(k,2m)$ are the Fourier coefficients of the Bernoulli polynomial $B_{2k}(t)$. The Fourier series of $B_{2k}(t)$, first computed by Hurwitz in 1890, is
\[
 B_{2k}(t) = \frac{(-1)^{k-1}(2k)!}{2^{2k-1} \pi^{2k}} 
 \sum_{m=1}^{\infty} \frac{\cos(2 \pi m t)}{m^{2k}},
 \qquad
 t \in [0,1),
 \ k \geq 1,
\]
and setting $t = 0$ yields \eqref{eq:sum2k}. However, this approach is not as simple or direct, since it necessitates the basic facts about Fourier series, not the least of which are the issue of pointwise convergence and the Riemann-Lebesgue Lemma. Fortunately, the latter is easily proved for $C^{1}$ functions via integration by parts (as we have done above). Pointwise convergence, as shown in~\cite{Cher}, can be dealt with using the ``telescoping trick''~\eqref{eq:telescope}, from which we can see the Fourier coefficients as differences of integrals involving the Dirichlet kernels $\sin((2n+1)\pi t)/\sin(\pi t)$.

The use of the same ``telescope'' \eqref{eq:telescope} in \cite{Ben} brought to our attention that this idea, used there to compute only $\zeta(2)$ via the integral $\int_{0}^{\pi} x \cos(mx) \, dx$, could also be used to compute $\zeta(2k)$. The crucial fact is realizing that, instead of powers, the natural functions to integrate against the cosine are the Bernoulli polynomials.

The integral for $\zeta(2k+1)$ in \eqref{eq:zeta-impar} is also known within the context of Fourier Analysis (see~\cite{Hauss}). It corresponds to the conjugate function of the Bernoulli polynomial $B_{2k+1}(t)$ evaluated at the origin.


\end{document}